\documentclass[12pt]{article}

\usepackage{amssymb}
\usepackage{amsmath}
\usepackage{theorem}

\newcommand{\Pf}{{\em Proof}. }
\newcommand{\EPf}{\hfill$\square$}
\newcommand{\copol}[2]{\mbox{copol}(#1,#2)}
\newcommand{\covar}[2]{\mbox{covar}(#1,#2)}
\newcommand{\D}{\mbox{$\mathcal D$}}
\newcommand{\E}{\mbox{$\mathcal E$}}
\newcommand{\inn}[2]{\mbox{$\mathcal{h} #1,#2 \mathcal{i}$}}
\newcommand{\vecsp}[1]{\mbox{$\mathcal{h} #1\mathcal{i}$}}
\newcommand{\Lg}{\mbox{$\mathfrak{g}$}}
\newcommand{\Lk}{\mbox{$\mathfrak{k}$}}
\newcommand{\Lh}{\mbox{$\mathfrak{h}$}}
\newcommand{\SU}[1]{\mbox{$\mathbf{SU}(#1)$}}
\newcommand{\U}[1]{\mbox{$\mathbf{U}(#1)$}}
\newcommand{\SP}[1]{\mbox{$\mathbf{Sp}(#1)$}}
\newcommand{\SO}[1]{\mbox{$\mathbf{SO}(#1)$}}

\newcommand{\Spin}[1]{\mbox{$\mathbf{Spin}(#1)$}}
\newcommand{\G}{\mbox{$\mathbf{G}_2$}}

\newtheorem{thm}{Theorem}

\newtheorem{rmk}[thm]{Remark}
\newtheorem{cor}[thm]{Corollary}
\newtheorem{prop}[thm]{Proposition}
\newtheorem{lem}[thm]{Lemma}

\begin{document}

\title{Copolarity of isometric actions}
\author{Claudio Gorodski\footnote{Partially supported
by CNPq grant 300720/93-9 and FAPESP grant 01/04793-8.}\hspace{.1cm} 
and Carlos Olmos\footnote{
Supported by Universidad Nacional de C\'ordoba and CONICET, partially 
supported by CIEM, Secyt-UNC and ANPCYT.}\hspace{.1cm}
 and Ruy Tojeiro\footnote{
Partially supported by CNPq grant~300229/92-5 and FAPESP grant
01/05318-1.}}


\maketitle

\begin{abstract}
We introduce a new integral invariant for isometric actions
of compact Lie groups, the \emph{copolarity}.
Roughly speaking, it measures how far from being polar
the action is. We generalize some results 
about polar actions in this context. 
In particular, we develop some of the structural
theory of copolarity $k$ representations,
we classify the irreducible representations of 
copolarity one,
and we relate the copolarity of an isometric action
to the concept of variational completeness
in the sense of Bott and Samelson.
\end{abstract}

\section{Introduction}

An isometric action of a compact Lie group $G$ on a 
complete Riemannian manifold $M$ is called \emph{polar}
if there exists a connected, complete submanifold
$\Sigma$ of $M$ which intersects all $G$-orbits 
and such that $\Sigma$ is orthogonal to every $G$-orbit it meets. 
Such a submanifold is called a \emph{section}. It is easy to see
that a section is automatically totally geodesic. If the 
section is also flat in the induced metric, then the action is called
\emph{hyperpolar}. 
In the case of Euclidean spaces, there is clearly no
difference between polar and hyperpolar representations since totally
geodesic submanifolds of an Euclidean space are affine
subspaces.
Polar representations were classified by Dadok~\cite{D}
and it follows from his work
that a polar representation of a compact 
Lie group is orbit equivalent to (i.~e.~has the same orbits as)
the isotropy representation of a symmetric space.

In this paper we introduce a new invariant for isometric actions
of compact Lie groups, the \emph{copolarity}.
Roughly speaking, it measures how far from being polar
the action is. This is based on the idea of 
a $k$-\emph{section}, which is a generalization of the concept of section.
The \emph{minimal $k$-section} passing through a regular point of the action
is the smallest connected, complete,
totally geodesic submanifold of the ambient space
passing through that point 
which intersects all the orbits and such that, 
at any intersection point with a principal orbit, its tangent space
contains the normal space of that orbit with codimension $k$.
It is easy to see that this is a good definition and uniquely specifies
an integer~$k$ which we call the \emph{copolarity}
of the isometric action (see Section~\ref{sec:k-sec}). 
It is also obvious that the $k=0$ case precisely
corresponds to the polar actions. 

It is apparent that for most isometric actions the minimal
$k$-section coincides with the ambient space. Note that in this case
$k$ equals the dimension of a principal orbit. 
We say that such isometric actions have \emph{trivial copolarity}.
The obvious questions that emerge are:
\begin{quote}
\em
What are the isometric actions with nontrivial copolarity?
What is the meaning of the integer $k$?
\end{quote}
In this paper we examine this problem in the case of orthogonal 
representations. 
Examples of representations of nontrivial copolarity
and minimal $k$-sections appear 
naturally in the framework of the reduction principle in
compact 
transformation groups (see~\cite{GS,SS,S1} for that principle). In fact,
in~\cite{GTh2} the reduction principle was used to describe 
the geometry of the irreducible representations 
in the table of Theorem~\ref{thm:1} below, which have 
copolarity $1$.
In that paper one was motivated by the
fact that the orbits of those representations
are tautly embedded in Euclidean space; we call representations
with this property \emph{taut}. 
This work is mainly motivated by the desire 
to better understand and generalize 
that description.

We give a complete answer to the above questions for the extremal 
values of the invariant $k$. 
Namely, let $(G,V)$ be an irreducible representation of a 
compact connected Lie group. Let $n$ be the dimension 
of a principal orbit. We prove the following two theorems. 

\begin{thm}\label{thm:1}
If $k=1<n$, then $(G,V)$ 
is one of the following orthogonal representations ($m\geq2$):
\[ \begin{array}{|c|c|}
\hline
\SO2\times\Spin9 & \mbox{(standard)}\otimes_{\mathbf R}\mbox{(spin)} \\
\U2\times\SP m & \mbox{(standard)}\otimes_{\mathbf C}\mbox{(standard)} \\
\SU2\times\SP m & \mbox{(standard)}^3\otimes_{\mathbf H}\mbox{(standard)} \\
\hline
\end{array}\]
\end{thm}

\begin{thm}\label{thm:2}
If $k=n-1$ or $k=n-2$, then $k=0$.
\end{thm}

Theorems~\ref{thm:1} and~\ref{thm:2} will be proved as 
corollaries of some other, stronger results,
see Corollaries~\ref{cor:one} and~\ref{cor:cod2} respectively. 
A couple of remarks are in order. 
Theorem~\ref{thm:2} says that in the nonpolar case a nontrivial
minimal $k$-section must have codimension at least $3$. 
Also, the three representations listed in the table of
Theorem~\ref{thm:1} are precisely the
irreducible representations of cohomogeneity~$3$ that are not 
polar (see~\cite{Y,D}). In fact, according to the main result
of~\cite{GTh3} (see also~\cite{GTh1}), these three
representations
together with the polar ones precisely comprise all the 
taut irreducible representations. Hence, we have the following
beautiful characterization of taut irreducible
representations.

\begin{thm}\label{thm:3}
An irreducible representation of 
a compact Lie group is taut if and only if $k=0$ or $k=1$. 
\end{thm}

Regarding Theorem~\ref{thm:3}, it is worth pointing out
that the case $k=1=n$ is impossible, for such a representation
would be orbit equivalent to a linear circle action 
and hence, by irreducibility, that would have to be
the standard action of $\SO2$ on $\mathbf R^2$,
which has $k=0$.
Notice that Theorems~\ref{thm:1} and~\ref{thm:3}
cease to hold if the representation 
is not irreducible 
as can be seen by taking 
the $7$-dimensional representation of $\U2$ given by the 
direct sum of the vector representation on $\mathbf C^2$
and the adjoint representation on $\mathfrak{su}(2)$
(It is interesting to remark that
this representation still has cohomogeneity $3$.)

Another result we would like to explain here is the following.
Let $(G,V)$ be an orthogonal
representation of nontrivial copolarity. 
It is easy to see that the $G$-translates of  
a nontrivial minimal $k$-section 
naturally determine a 
group invariant foliation $\mathcal F$
on the $G$-regular set of $V$
(in fact, here the $k$-section need not 
be minimal, but we do not go into details in this introduction). 
We prove the following theorem
(see Theorem~\ref{thm:int}).

\begin{thm} 
If the distribution orthogonal to $\mathcal F$ 
is integrable, 
then $(G,V)$ is orbit equivalent to a
direct product representation $(G_1\times G_2,V_1\oplus V_2)$, where 
$G_1$, $G_2$ are subgroups of $G$, $(G_1,V_1)$ is a polar
representation and $(G_2,V_2)$ is any representation;
here the leaves of the distribution 
orthogonal to $\mathcal F$ correspond to the $G_1$-orbits.
In particular, if $(G,V)$ is nonpolar then it cannot be irreducible. 
\end{thm}

In this paper we also relate the copolarity of an isometric action 
to the concept of variationally complete actions which was
introduced by Bott in \cite{B} (see also~\cite{BS}).
Roughly speaking, an isometric action of a compact 
Lie group on a complete Riemannian manifold is
\emph{variationally complete} if it produces enough Jacobi fields along 
geodesics to determine the multiplicities of focal points to the orbits
(see Section~\ref{sec:varcomp} for the precise definition).  
Conlon proved in~\cite{C} that a hyperpolar action of a compact Lie
group on a complete Riemannian manifold is variationally complete. 
On the other hand, it is known that 
a variationally complete representation
is polar~\cite{DO,GTh1}, and that a variationally complete action on 
a compact symmetric space is hyperpolar~\cite{GTh4}. 
This implies that
the converse to Conlon's theorem is true for actions on Euclidean 
spaces or compact symmetric spaces. 
In this paper we introduce the notion of variational
co-completeness of an isometric action and prove that
it does not exceed $k$ for an action that admits a flat $k$-section
(Theorem~\ref{thm:Conlon}). This reduces to Conlon's 
theorem for $k=0$. We also prove a weak converse of this result
in the case of representations (Theorem~\ref{thm:converse}). 

The paper is organized as follows. 
We first define $k$-sections 
and the copolarity of an isometric action (Section~\ref{sec:k-sec})
and present some examples (Section~\ref{sec:ex}). 
Then we go on to develop some of the structural theory
of copolarity $k$ representations. In particular, we
show that the copolarity of an orthogonal representation 
behaves well with respect to taking slice representations
(Theorem~\ref{thm:slice-copolar})
and forming direct sums (Theorem~\ref{thm:red}), 
and we obtain a reduction principle in terms of 
$k$-sections (Theorem~\ref{thm:reduction}). 
We also show that the codimension of a nontrivial
minimal $k$-section of an irreducible representation
is at least $3$ (Corollary~\ref{cor:cod2}),
and characterize the orthogonal representations 
admitting a minimal $k$-section whose orthogonal 
distribution is integrable (Theorem~\ref{thm:int}). 
We describe the geometry of a principal orbit
of a representation of copolarity one (Theorem~\ref{thm:one})
and classify the irreducible representations of copolarity one
(Corollary~\ref{cor:one}). 
We finally prove the extension 
of Conlon's theorem for copolarity $k$ actions (Theorem~\ref{thm:Conlon})
and its weak converse in the case of representations
(Theorem~\ref{thm:converse}). As a corollary, we generalize
a result about polar representations (Corollary~\ref{cor:normal}).

As a final note, we recall that 
the principal orbits of polar representations
can be characterized as being the only compact homogeneous 
isoparametric submanifolds of Euclidean space~\cite{PT}. 
An open problem in the area is to similarly 
characterize the principal orbits of 
more general orthogonal representations in terms of their submanifold 
geometry and topology. 
We believe that orthogonal representations of low copolarity 
may serve as testing cases for this problem. 

The first author wishes to thank Prof.~Gudlaugur Thorbergsson
for very useful conversations. 
Part of this work was completed while the
third author was visiting
University of S\~ao Paulo (USP), 
for which he wishes to thank Prof.~Ant\^onio
Carlos Asperti and the other colleagues from USP
for their hospitality, and FAPESP for financial support.

\section{Actions admitting $k$-sections}\label{sec:k-sec}
\setcounter{thm}{0}

Let $(G,M)$ be an isometric action of the compact Lie group $G$ on the
complete Riemannian manifold $M$. A $k$-\emph{section} for 
$(G,M)$, where $k$ is a nonnegative integer,
is a connected, complete submanifold $\Sigma$ of $M$ such that
the following hold:
\begin{enumerate}
\item[(C1)] $\Sigma$ is totally geodesic in $M$;
\item[(C2)] $\Sigma$ intersects all $G$-orbits;
\item[(C3)] for every $G$-regular point $p\in\Sigma$ we have that
 $T_p\Sigma$ contains the normal space $\nu_p(Gp)$
as a subspace of codimension $k$;
\item[(C4)] for every $G$-regular point $p\in\Sigma$ we have that
if $gp\in\Sigma$
for some $g\in G$, then $g\Sigma=\Sigma$.
\end{enumerate}
If $\Sigma$ is a $k$-section through $p$, 
then $g\Sigma$ is a 
$k$-section through $gp$ for any $g\in G$.
We also remark that:
since a $k$-section $\Sigma$ is connected, complete
and totally geodesic, for every $p\in\Sigma$ we have that
$\Sigma=\exp_p T_p\Sigma$; and, since the $G$-orbits are compact,
for every $p\in M$ we have that the set $\exp_p\nu(Gp)$ intersects
all $G$-orbits. Using these remarks, it
is easy to see that, given a $G$-regular $p\in M$,
the connected component containing $p$ of the
intersection of a $k_1$-section 
and a $k_2$-section passing through $p$ is a
smooth submanifold and it 
is a $k$-section passing through $p$ with $k\leq\min\{k_1,k_2\}$.
It is also clear that the ambient $M$ is a trivial $k$-section
for $k$ equal to the dimension of a principal orbit. 
It follows that the set of $k$-sections, $k=0,1,2,\ldots$, passing 
through a fixed regular point admits a unique minimal element.
We say that the \emph{copolarity} of $(G,M)$ is $k_0$,
and we write $\copol{G}{M}=k_0$, if that minimal element
is a $k_0$-section. 
In this way, the copolarity is well defined for 
any isometric action $(G,M)$ as being an integer~$k_0$ between 
zero and the dimension of a principal orbit,
and then a $k_0$-section is uniquely determined
through any given regular point. We say that $(G,M)$ has 
\emph{nontrivial copolarity} if $k_0$ is 
strictly less than the dimension of a principal orbit
(or, equivalently, if a $k_0$-section is properly
contained in $M$). 
Since a $0$-section is simply a section, we have that
$\copol{G}{M}=0$ if and only if $(G,M)$ is polar.
Note that the set of connected, complete submanifolds of $M$ passing 
through a fixed $G$-regular point and satisfying 
only conditions (C1), (C2) and (C3) is also closed under connected 
intersection, and that a minimal element in this set
automatically satisfies condition (C4), so that it represents the same
minimal $k_0$-section. \emph{It follows from this observation that,
in order to show that an isometric action has copolarity at most $k$,
it is enough to construct a connected, complete submanifold of codimension
$k$ satisfying conditions (C1), (C2) and (C3).}
On the other hand, note that many of our applications will not depend 
on the fact that a $k$-section is minimal, but rather will
depend on the fact that it satisfies condition (C4). 

Next we discuss the conditions in the
definition of a $k$-section. Note that
if $k=0$ condition~(C4) is unnecessary and also~(C1) follows
rather easily from~(C2) and~(C3), but in general we cannot 
dispense with them. A standard argument also shows that
in the general case condition~(C2) follows from~(C3), if the latter is 
not empty, namely if we assume that $\Sigma$ contains a regular point. 
Condition~(C4) (combined with~(C2))
is equivalent to the fact that the $G$-translates of $\Sigma$
define a foliation in the regular set of $M$.  
Even more interesting is the following rephrasing of condition~(C4).
(Note that condition~(C3) implies that the intersection of a $k$-section 
with a principal orbit is a smooth manifold.)

\begin{prop}\label{prop:distr}
Let $(G,M)$ be an isometric action of the compact Lie group $G$ on the
complete Riemannian manifold $M$.
Suppose $\Sigma$ be a connected, complete submanifold of $M$ 
satisfying conditions~(C1), (C2) and~(C3). For every regular 
$p\in\Sigma$ define a $k$-dimensional subspace $\D_p$ of $T_p(Gp)$ by
$\D_p=T_p\Sigma\cap T_p(Gp)$.  
Then condition~(C4) is equivalent to any one of the following assertions:
\begin{enumerate}
\item[(a)] the subspaces $\{\D_p\}$ extend
to a $k$-dimensional, $G$-invariant 
distribution $\D$ on the regular set of $M$;
\item[(b)] for every $G$-regular point $p\in\Sigma$ 
we have that if $gp\in\Sigma$ 
for some $g\in G$, then $g_*\D_p=\D_{gp}$;
\item[(c)] there exists a principal orbit $\xi$ 
such that if $p$, $gp\in\Sigma\cap\xi$ for some
$g\in G$, then $g_*\D_p=\D_{gp}$.
\end{enumerate}
Moreover, any one of the above conditions implies that
the stabilizer $G_\Sigma$ of $\Sigma$ 
acts transitively on the intersection of 
$\Sigma$ with any principal orbit.
\end{prop}

\Pf Since $\Sigma$ is connected, complete and totally geodesic in $M$, 
for any $G$-regular point $p\in\Sigma$ we have that
$\Sigma=\exp_p T_p\Sigma$, where $T_p\Sigma=\D_p\oplus\nu_p(Gp)$.
We now show that~(c) implies~(C4).
Assume that~(c) holds with 
respect to the principal orbit $\xi$. Let $q\in\Sigma$ 
be a $G$-regular point. There exists a minimal geodesic 
from~$q$ to~$\xi$. Therefore we can write
$p=\exp_q v$ for some $p\in\xi$ and $v\in\nu_q(Gq)$. 
Conditions~(C3) and~(C1) imply that $p\in\Sigma$. 
Now if $gq\in\Sigma$ for some $g\in G$, then 
$gp=\exp_{gq}g_*v\in\exp_{gq}(\nu_{gq}(Gq))\subset
\exp_{gq}(T_{gq}\Sigma)=\Sigma$. 
Since $p$, $gp\in\Sigma\cap\xi$, by~(c) we deduce that $g_*\D_p=\D_{gp}$.
It follows that $g\Sigma=g\exp_p(\D_p\oplus\nu_p(Gp))=
\exp_{gp}(g_*\D_p\oplus g_*\nu_p(Gp))=
\exp_{gp}(\D_{gp}\oplus\nu_{gp}(Gp))=\Sigma$, and this 
gives~(C4). 

Notice that~(a) is just a reformulation of~(b), and
the fact that $G_\Sigma$ is transitive on the intersection of 
$\Sigma$ with any principal orbit follows immediately from~(C4).
Also, the equivalence between~(C4) and~(b) 
follows from $\Sigma=\exp_p(\D_p\oplus\nu_p(Gp))$.
Since~(b) trivially implies~(c),
this completes the proof of the proposition. \EPf

\medskip

Let $\Sigma$ be a $k$-section of $(G,M)$ and 
consider the distribution~$\D$ as in  Proposition~\ref{prop:distr}.
Define $\E$ to be the distribution on the regular 
set of $M$ such that $\E_p$ is the orthogonal complement 
of $\D_p$ in $T_p(Gp)$. If $p$ 
is $G$-regular, then it follows from the fact that 
$\Sigma$ is totally geodesic in $M$ that 
the distribution $\D|_{Gp}$ is autoparallel
in $Gp$ and invariant under the 
second fundamental form $\alpha$ of $Gp$, in the sense
that $\alpha(\D,\E)=0$ on $Gp$. 

The following proposition will be used later 
to show that some orthogonal 
representations admit $k$-sections.

\begin{prop}\label{prop:suf}
Let $(G,V)$ be an orthogonal representation of a compact Lie 
group on a Euclidean space $V$. Let $Gp$ be a principal
orbit. Suppose there is a $G$-invariant, autoparallel,
$k$-dimensional distribution $\D$ on $Gp$ such that $\alpha(\D,\E)=0$,
where $\alpha$ is the second fundamental form of $Gp$ and
$\E$ is the distribution on $Gp$ orthogonal to $\D$. 
Define $\Sigma=\D_p\oplus\nu_p(Gp)$.
\begin{enumerate}
\item[(a)] Suppose that for every $v\in\nu_p(Gp)$ with 
$q=p+v$ a $G$-regular point we have that $\E_p\subset T_q(Gq)$.
Then $\Sigma$ satisfies conditions~(C1), (C2) and~(C3)
in the definition of a $k$-section. 
\item[(b)] Suppose that, in addition to the hypothesis in~(a),
for every $v\in\nu_p(Gp)$ with 
$q=p+v\in Gp$ we have that $\E_p=\E_q$.
Then $\Sigma$ is a $k$-section for $(G,V)$. 
\end{enumerate}
\end{prop}

\Pf Let $\beta$ be a the maximal integral submanifold of $\D$ through
$p$. Since $\beta$ is totally geodesic in $Gp$ and $\alpha(\D,\E)=0$,
we have that the covariant derivative in $V$ of a $\E$-section
along $\beta$ is in $\E$, which implies that
$\E$ is constant in $V$ along $\beta$. 
Therefore $\beta$ is contained in the affine subspace 
of $V$ orthogonal to $\E_p$ which is $\Sigma$ and, in fact,
$\beta$ is the connected component of $\Sigma\cap Gp$ 
containing $p$. 

Now if $gp\in\beta$ then $\E_p=\E_{gp}$. Taking 
orthogonal complements in $V$, we get that 
$\Sigma=\D_{gp}\oplus\nu_{gp}(Gp)$; 
but the right hand side in turn equals
$g\Sigma$, as $\D$ is $G$-invariant.
This shows that conditions~(C3) and~(C4) are already verified
for points in~$\beta$.

Next let $\gamma\neq\beta$ be a connected component of the intersection 
of $\Sigma$ with a principal orbit (which possibly could be $Gp$). 
Let $q\in\gamma$ and consider the minimal geodesic in $V$ from $q$
to $\beta$. Then we can write $q=gp+v$, where $gp\in\beta$ for some
$g\in G$ and $v\in\nu_{gp}(Gp)$. Since $g^{-1}q=p+g^{-1}v$, 
we have $T_q(Gq)=gT_{g^{-1}q}(Gq)\supset
g\E_p=\E_{gp}=\E_p$ (where the inclusion follows from the hypothesis
in~(a)).  
Taking orthogonal complements in $V$, we get that 
$\nu_q(Gq)\subset\Sigma$. 
This shows that condition~(C3) is fully verified. 
If, in addition, $\gamma$ is a connected 
component of $\Sigma\cap Gp$ and $q=hp$ for some $h\in G$,
then we can write $h\E_p=\E_q=g\E_{g^{-1}q}=
g\E_p=\E_{gp}=\E_p$ (where the equality $\E_{g^{-1}q}=\E_p$
follows from the hypothesis in~(b)).  
Taking orthogonal complements we get $h\Sigma=\Sigma$ which,
by Proposition~\ref{prop:distr}(c), finally implies 
condition~(C4). \EPf

\section{Examples}\label{sec:ex}
\setcounter{thm}{0}

In this section we exhibit some examples of actions admitting 
$k$-sections.

\subsection{Product actions}

Let $(G_1,M_1)$ be a polar action with a section $\Sigma_1$,
and let $(G_2,M_2)$ be any isometric 
action. Let $G=G_1\times G_2$, $M=M_1\times M_2$ and consider the
product action $(G,M)$. Then it is immediate to see
that $\Sigma=\Sigma_1\times M_2$ is a 
$k$-section of $(G,M)$, where $k$ is the dimension of a principal
orbit of $(G_2,M_2)$. Note that in this example the distribution 
$\E$ is integrable and its leaves are the $G_2$-orbits in $M_2$.
We will show later that, in the case of orthogonal
representations, this example is essentially the only one 
whose distribution $\E$ is integrable. 
 
\subsection{The reduction principle in compact
transformation groups}

Let $(G,M)$ be an isometric action and fix a principal isotropy
subgroup $H$. Then the $H$-fixed point set $\Sigma=M^H$ is a
$k$-section for $(G,M)$, where $k$ is the difference between the 
dimension of $\Sigma$ and the cohomogeneity of $(G,M)$. 
In fact, $\Sigma$ is a totally geodesic submanifold of $M$, being
the common fixed point set of a set of isometries of $M$. 
This is condition (C1) in the definition of a $k$-section. 
Condition~(C3) follows from the fact that the slice representation 
at a regular point is trivial. Moreover, if $p$, $gp\in\Sigma$ are regular 
points for some $g\in G$, then both the isotropy subgroups at $p$ and $gp$ 
are $H$, so $g$ normalizes $H$ and therefore fixes $\Sigma$. 
This is condition~(C4). The rest follows. 

It is interesting to notice that sometimes the group $G$ 
can be enlarged to another group $\hat G$ 
that has the same orbits as $G$ on $M$
but has a larger principal isotropy subgroup $\hat H$, and
then the $\hat H$-fixed point set is smaller that the 
$H$-fixed point set (see~\cite{S1}). 

This example is very important in the sense 
it shows that if $(G,M)$ is an arbitrary isometric action,
and $\Sigma$ is $k$-section which is \emph{minimal}
(so that $\copol{G}{M}=k$),
then we always have that $\Sigma\subset M^{G_p}$,
for any $G$-regular point $p\in\Sigma$. Therefore $G_p$ 
acts trivially on $\D_p=T_p\Sigma\cap T_p(Gp)$.
We do not know of any example of an isometric action
such that the minimal $k$-section is \emph{strictly} contained in the 
fixed point set of a principal isotropy subgroup. 

In any case, $G_p\subset G_\Sigma$ for a minimal $k$-section
$\Sigma$ and a $G$-regular point $p\in\Sigma$. 
This implies that the isotropy subgroup of $G_\Sigma$ 
at~$p$ is $G_p$. Since $G_p$ acts trivially on $\Sigma$,
using Lemma~\ref{lem:inj} below, we get that 
$Gp\cap\Sigma=G_\Sigma p=G_\Sigma/G_p$ is a group 
manifold.  

\subsection{Some orthogonal representations of low 
copolarity}\label{subsec:low}

Some calculations in~\cite{GTh3,GTh1} involving the reduction principle in 
transformation groups produced some examples of 
irreducible orthogonal representations of low 
copolarity. 
\begin{enumerate}
\item[(i)] The three irreducible representations of cohomogeneity $3$
have copolarity $1$. These are listed in the table of
Theorem~\ref{thm:1}.
\item[(ii)] The tensor product of the
vector representation of $\SO3$ and the $7$-dimensional representation of 
$\G$ is an irreducible
representation of copolarity $2$. 
\item[(iii)] Let $(S^1\times H, V)$ be a polar irreducible representation
which is \emph{Hermitian}, namely leaves a complex structure on $V$ 
invariant. Assume that the restricted representation $(H,V)$
is not orbit equivalent to $(S^1\times H, V)$.
Then $(H,V)$ has nontrivial copolarity $k>0$. 
(There are four families of such Hermitian
polar irreducible representations, and
the simplest example is maybe the action of the unitary group 
$\U3$ on the space of complex symmetric bilinear forms in three
variables. The induced representation
of $\SU3$ has real dimension $12$, cohomogeneity $4$
and copolarity $2$.)
\end{enumerate}

\subsection{Examples derived from polar actions}

Let $(G,V)$ be a polar representation of a compact 
Lie group $G$ on a Euclidean space. Then, of course, 
$\copol{G}{M}=0$. Nevertheless, in general there are 
interesting examples of $k$-sections for $(G,V)$ with $k>0$. 

In fact, let $Gp$ be a principal orbit. It is known that 
equivariant normal vector fields along $Gp$ are parallel
in the normal connection and therefore
$Gp$ is an \emph{isoparametric submanifold} of $V$, namely 
the principal curvatures of $Gp$ along a parallel normal 
vector field are constant and the normal bundle of $Gp$ in $V$
is globally flat. 
Let $v\in\nu_p(Gp)$, $A_{v}$ the Weingarten operator
in the direction of $v$ and 
$\lambda$ a nonzero principal curvature in the direction
of $v$ with multiplicity $k$. 
Since the subspace $\ker(A_v-\lambda\mbox{id}_{T_p(Gp)})$
of $T_p(Gp)$ is $G_p$- and $A_v$-invariant, it extends to a
$G$-invariant distribution $\D$ on $Gp$ which is
invariant under the second fundamental form of $Gp$, and one can
see from the Codazzi equation 
that $\D$ is also autoparallel. Moreover, 
it is very easy to deduce from the theory of isoparametric 
submanifolds
that the conditions on $\E$ from Proposition~\ref{prop:suf}
are satisfied, so $\Sigma=\D_p\oplus\nu_p(Gp)$ is a $k$-section.

\section{Structural theory of actions admitting $k$-sections}
\setcounter{thm}{0}

\subsection{Slice representations}

In this section we will prove that the copolarity of a
slice representation is not bigger than the copolarity of the 
original representation (Theorem~\ref{thm:slice-copolar}). 
Let $(G,M)$ be an isometric action
of the compact Lie group $G$ on the complete
Riemannian manifold $M$. Let $q\in M$. Then $G_q$ acts on
$T_qM$ via the differential, and $T_qM=T_q(Gq)\oplus\nu_q(Gq)$ is an
invariant decomposition. The orthogonal representation
$(G_q,\nu_q(Gq))$ is called the \emph{slice representation at~$q$}. 

\begin{lem}\label{lem:slice-reg}
The following assertions are equivalent:
\begin{enumerate}
\item[(a)] $v\in\nu_q(Gq)$ is $G_q$-regular.
\item[(b)] There exists $\epsilon>0$ such that
$\exp_q(tv)$ is $G$-regular for $0<t<\epsilon$.
\item[(c)] $\exp_q(t_0v)$ is $G$-regular for some $t_0>0$. 
\end{enumerate}
\end{lem}

\Pf (c) implies (b): Let $p=\exp_q(t_0v)$ be $G$-regular. Then
there exists $\epsilon>0$ such that 
\[ G_{\exp_q(tv)}=(G_q)_{tv}=(G_q)_v\subset G_p\subset G_q, \]
for $0<t<\epsilon$, where the first equality follows from the fact that
the exponential map is an equivariant diffeomorphism of a
small normal disk of radius $\epsilon$ onto the image, 
and the last inclusion follows from the fact that 
the slice representation at $p$ is trivial as $p$ is $G$-regular.  
Again by the $G$-regularity of~$p$, we have that
$G_{\exp_q(tv)}=G_p$ for $0<t<\epsilon$,
and hence $\exp_q(tv)$ is $G$-regular for $0<t<\epsilon$.

Clearly (b) implies (c), and the equivalence of (a) with (b) 
follows from the equality $G_{\exp_q(tv)}=(G_q)_v$ for $0<t<\epsilon$
and the slice theorem. \EPf

\medskip

In the following we consider the case of an orthogonal 
representation $(G,V)$. 
Let $q\in V$, choose a $k$-section $\Sigma$ through $q$ and
consider the slice representation $(G_q,\nu_q(Gq))$. 

\begin{lem}
There exists $v\in\Sigma\cap\nu_q(Gq)$ which is a $G_q$-regular point. 
\end{lem}

\Pf Let $\xi$ be a principal $G$-orbit and choose a connected component 
$\beta$ of $\Sigma\cap\xi$. Let $c(t)=q+tv$, $0\leq t\leq l$,
be a minimal geodesic in $\Sigma$ from 
$q$ to $\beta$. Then 
$\dot{c}(l)\in\Sigma=T_{c(l)}\beta\oplus\nu_{c(l)}(Gc(l))$.
But $\dot{c}(l)$ must be orthogonal to $\beta$, by minimality.
Therefore $\dot{c}(l)\in\nu_{c(l)}(Gc(l))$. Since a geodesic orthogonal
to an orbit must be orthogonal to every orbit it meets, 
$v=\dot{c}(0)\in\Sigma\cap\nu_q(Gq)$.  
It follows from Lemma~\ref{lem:slice-reg}
that $v$ is $G_q$-regular. \EPf

\medskip

For the $G$-regular $p=q+tv\in V$, $0<t<\epsilon$, we then have that 
$\Sigma=\nu_p(Gp)\oplus\D_p$, where $\D$ has rank $k$. 

\begin{lem}
There is an orthogonal decomposition 
\[ T_p(G_q p)= T_p(G_q p) \cap\D_p\oplus T_p(G_q p)\cap\E_p. \]
\end{lem}

\Pf Note that 
$T_p(Gp)=\D_p\oplus\E_p$ is an $A_v$-invariant decomposition,
where $A_v$ denotes the Weingarten operator of $Gp$ with respect
to the, say, unit vector $v\in\nu_p(Gp)$.  
We first claim that $T_p(G_q p)$ is contained in the eigenspace of $A_v$ 
corresponding to the eigenvalue $-1/t$. In fact,
let $u=\frac{d}{ds}|_{s=0}\varphi_s(p)\in T_p(G_q p)$ for 
$\varphi_s\in G_q$, $\varphi_0=1$. Define $\hat{v}(s)=\varphi_sv$ 
normal vector field along $s\mapsto\varphi_sp$. Then 
$\hat v(s)=\frac{1}{t}(\varphi_sp-q)$ and 
\[ -A_vu+\nabla_u^\perp\hat{v}=\frac{d}{ds}
\Big|_{s=0}\hat{v}(s)=\frac{1}{t}u, \]
so, by taking tangent components, we get that $A_vu=-\frac{1}{t}u$. 

Next write $u=u'+u''$ where $u'\in\D_p$ and $u''\in\E_p$. 
Then $u'$ and $u''$ are eigenvectors of $A_v$ with eigenvalue $-1/t$. 
Since $u''\in\E_p$, we have that $u''\in T_p(G_q p)$
by Corollary~\ref{cor:partial-vc}. 
Therefore, we also have that $u'\in T_p(G_q p)$. \EPf

\medskip

Let $\D_{1p}=T_pG_q(p) \cap\D_p$ and 
$\E_{1p}=T_pG_q(p) \cap\E_p$.
Let $\D_{2p}$ be the orthogonal complement of $\D_{1p}$ in $\D_p$
and let $\E_{2p}$ be the orthogonal complement of $\E_{1p}$ in 
$\E_p$. 

\begin{lem}\label{lem:d2perp}
We have that $\E_{2p}\subset T_q(Gq)$.
\end{lem}

\Pf Let $u=\frac{d}{ds}|_{s=0}\varphi_sp\in\E_{2p}$ with 
$\varphi_s\in G$, $\varphi_0=1$. Without loss of generality we may assume
that $A_vu=\lambda u$, and then $\lambda\neq-1/t$ because the $-1/t$
eigenvectors in $\E_p$ belong to $T_p(G_q p)$ 
(Corollary~\ref{cor:partial-vc}).
Let $r:Gp\to Gq$ be the canonical equivariant
submersion. We compute
\[ r_*(u)=\frac{d}{ds}\Big|_{s=0}\underbrace{r\varphi_s}_{=\varphi_sr} (p)  
       =\frac{d}{ds}\Big|_{s=0}\underbrace{\varphi_s (q)}_{=\varphi_s(p)
-t\varphi_s(v)} 
       =u-t\frac{d}{ds}\Big|_{s=0}\hat v(s), \]
so 
 \[ u=r_*(u)+t(-A_vu+\nabla_u^\perp\hat v). \]
Now $\nabla_u^\perp\hat v=0$ because $u\in\E_p$ 
(Corollary~\ref{cor:partial-parallel}). 
Hence $u=(1+\lambda t)^{-1}r_*(u)\in T_q(Gq)$. \EPf 

\medskip

Notice that each one of $\D_{1p}$, $\E_{1p}$, $\nu_p(Gp)\oplus\D_{2p}$ 
and $\E_{2p}$ 
is a constant subspace of $V$ with respect to $t\in(0,\epsilon)$, 
where~$p=q+tv$. 

\begin{lem}\label{lem:slice-intersect}
We have that $\Sigma\cap\nu_q(Gq)$ intersects all $G_q$-orbits.
\end{lem}

\Pf It follows from Lemma~\ref{lem:d2perp} that 
$\nu_q(Gq)\subset\Sigma\oplus\E_{1p}$. 
Let $N_v$ be the normal space to $T_p(G_q p)$ in $\nu_q(Gq)$. 
Then $N_v$ is orthogonal to $\E_{2p}$ so that
\[ \nu_q(Gq) = \underbrace{T_p(G_q p)}_{=\D_{1p}\oplus\E_{1p}} 
 \oplus\underbrace{N_v}_{\subset\nu_p(Gp)\oplus\D_{2p}}. \]
Hence $\nu_q(Gq)\cap\Sigma=N_v\oplus\D_{1p}$. Since $\nu_q(Gq)\cap\Sigma$ 
contains $N_v$ and $v$ is $G_q$-regular, we get that 
$\nu_q(Gq)\cap\Sigma$ intersects all $G_q$-orbits. \EPf

\medskip

For a $G_q$-regular $w\in\nu_q(Gq)\cap\Sigma$, 
we have that $q+tw$ is $G$-regular for $t>0$ small,
and thus it makes sense to define 
$\D_{1w}=\D_{q+tw}\cap T_{q+tw}G_q(q+tw)$.
It is clear that $\D_1$ is a $G_q$-invariant distribution 
on the $G_q$-regular set of $\nu_q(Gq)$. 
Moreover, $\nu_q(Gq)\cap\Sigma=N_w\oplus\D_{1w}$, as above.
It follows that $\nu_q(Gq)\cap\Sigma$ is a $k_1$-section
for $(G_q,\nu_q(Gq))$, where $k_1$ is the dimension of $\D_1$. 
Since the dimension of $\D_1$ is not bigger than the dimension
of $\D$, we finally conclude:

\begin{thm}\label{thm:slice-copolar}
If $\mbox{copol}(G,V)\leq k$, then 
$\mbox{copol}(G_q,\nu_q(Gq))\leq k$.
\end{thm}

\subsection{The reduction}

In this section we establish a reduction principle
for orthogonal representations in terms of $k$-sections
(Theorem~\ref{thm:reduction}).
Let $(G,V)$ be an orthogonal representation
admitting a $k$-section $\Sigma$.

\begin{lem}\label{lem:trans-sec}
Let $q\in V$. Then the isotropy subgroup $G_q$ is transitive on the set
$\mathcal F$ of $k$-sections through $q$ which are 
$G$-translates of $\Sigma$, 
namely $\mathcal F=\{g\Sigma:g\in G,\, q\in g\Sigma\}$.
\end{lem}

\Pf The result is trivial for a $G$-regular point $q$, since in this 
case there is a unique element in $\mathcal F$ by condition (C4). 
Assume $q$ is 
not a $G$-regular point, and let $\Sigma_1$, $\Sigma_2\in\mathcal F$. 
Let $S$ be a normal slice 
at $q$ and choose a $G_q$-regular $p\in S$. 
By Lemma~\ref{lem:slice-intersect}, $\Sigma_i\cap\nu_q(Gq)$ 
intersects all $G_q$-orbits on $S$, 
$i=1$, $2$. Therefore we can select 
$h_i\in G_q$ such that $h_ip\in\Sigma_i$. 
Now~$p\in h_1^{-1}\Sigma_1\cap h_2^{-1}\Sigma_2$ where
$p$ is $G$-regular, so
$h_2^{-1}h_1\Sigma_1=\Sigma_2$ where $h_2^{-1}h_1\in G_q$. \EPf

\begin{lem}\label{lem:inj}
Let $G_\Sigma$ be the stabilizer of $\Sigma$. Then 
$G_\Sigma q=Gq\cap\Sigma$, for all $q\in\Sigma$.
\end{lem}

\Pf Let $p$, $q\in\Sigma$ be in the same $G$-orbit. 
We need to show that they are in the same $G_\Sigma$-orbit. 
If they belong to a principal $G$-orbit, then the 
result follows from Proposition~\ref{prop:distr}. Suppose $p$ and $q$ are not 
$G$-regular points and $q=gp$ for some $g\in G$. 
Then $q\in\Sigma\cap g\Sigma$.
By Lemma~\ref{lem:trans-sec}, there is $h\in G_q$ such that 
$h\Sigma=g\Sigma$. Then $h^{-1}g\in G_\Sigma$ and $(h^{-1}g)p=q$. \EPf 

\begin{thm}\label{thm:reduction}
The inclusion $\Sigma\to V$ induces a homeomorphism 
$\Sigma/G_\Sigma\to V/G$.
\end{thm}

\Pf The surjectivity is implied by the fact that $\Sigma$ 
intersects \emph{all} $G$-orbits. The injectivity
is precisely Lemma~\ref{lem:inj}. Since the inclusion
$\Sigma\to V$ is continuous, the induced map 
$\Sigma/G_\Sigma\to V/G$ is continuous. 
The restriction $S(\Sigma)/G_\Sigma\to S(V)/G$,
where $S(\Sigma)\subset\Sigma$, $S(V)\subset V$ are the unit spheres,
is also a continuous bijection, and thus a homeomorphism,
as $S(\Sigma)/G_\Sigma$ is compact. Since $\Sigma/G_\Sigma\to V/G$
is the cone over $S(\Sigma)/G_\Sigma\to S(V)/G$, it follows that 
$\Sigma/G_\Sigma\to V/G$ is a homeomorphism.
\EPf

\medskip
In the remaining of this section we study the intersection of 
the $k$-section with a nonprincipal orbit. Let $q\in\Sigma$
be a singular point. 

\begin{lem}\label{lem:sigma-decomp}
We have an orthogonal decomposition
$\Sigma=\Sigma\cap T_q(Gq)\oplus\Sigma\cap\nu_q(Gq)$.
\end{lem}

\Pf Write
$\Sigma=T_q(G_\Sigma q)\oplus\nu^\Sigma_q(G_\Sigma q)$, where
$\nu^\Sigma_q(G_\Sigma q)$ denotes the normal space to 
$G_\Sigma q$ at $q$ in $\Sigma$. Lemma~\ref{lem:inj}
implies that $T_q(G_\Sigma q)=\Sigma\cap T_q(Gq)$. 
We still need to show that $\nu^\Sigma_q(G_\Sigma q)=\Sigma\cap\nu_q(Gq)$
in order to complete the proof. 

Since $\Sigma\cap\nu_q(Gq)\subset\nu^\Sigma_q(G_\Sigma q)$ and 
$\Sigma\cap\nu_q(Gq)$ intersects all $(G_q,\nu_q(Gq))$-orbits
by Lemma~\ref{lem:slice-intersect}, we have that there is 
$v\in\nu^\Sigma_q(G_\Sigma q)$ which is $G_q$-regular.
Now by Lemma~\ref{lem:slice-reg} we know that 
$p=q+t_0v$ is $G$-regular for some $t_0>0$.
We use the remark that a geodesic orthogonal to an orbit is 
orthogonal to every orbit it meets; we must have  
$v\in\nu_p^\Sigma(G_\Sigma p)$. Since $p$ is $G$-regular, 
this implies that $v\in\nu_p(Gp)$. Again, by the same remark,
$v\in\nu_q(Gq)$. \EPf

\begin{lem}\label{lem:sigma-sing}
Let $r:Gp\to Gq$ be the canonical equivariant
submersion, where $p$ is a $G$-regular point in the normal slice at $q$. 
Then $\Sigma\cap T_q(Gq)=r_*\D_{2p}$ and this is 
the orthogonal complement of $r_*\E_{2p}=\E_{2p}$ in 
$T_q(Gq)$.
\end{lem}

\Pf Write $p=q+tv$ where $v\in\nu_q(Gq)$.
Let $u\in\D_{2p}\subset T_p(Gp)$. A computation similar to the 
one done in Lemma~\ref{lem:d2perp} shows that
\[ r_*(u)=\underbrace{(\mbox{id}+tA_v)u}_{\in\D_2p}
+\underbrace{t\nabla^\perp_u\hat v}_{\in\nu_p(Gp)} \]
is a vector in $\Sigma$ and therefore orthogonal to $\E_{2p}$. 
Now 
\[ T_q(Gq)=r_*T_p(Gp)=r_*\D_p+r_*\E_p=
\underbrace{r_*\D_{2p}}_{\subset\Sigma}+\underbrace{r_*\E_{2p}}_{=\E_{2p}}. \]
Since $\Sigma$ is orthogonal to $\E_{2p}$, we conclude that
the last sum is an orthogonal direct sum. \EPf

\begin{cor}
We have that $r_*\D_{2p}\oplus N_v=\D_{2p}\oplus\nu_p(Gp)$.
\end{cor}

\Pf Consider the orthogonal decomposition
$\Sigma=\Sigma\cap T_q(Gq)\oplus\Sigma\cap\nu_q(Gq)$
from Lemma~\ref{lem:sigma-decomp}. On the one hand, we know that
$\Sigma=\nu_p(Gp)\oplus\D_{2p}\oplus\D_{1p}$.
On the other hand, we have that $\Sigma\cap T_q(Gq)=r_*\D_{2p}$
by Lemma~\ref{lem:sigma-sing} and $\Sigma\cap\nu_q(Gq)=N_v\oplus\D_{1p}$ 
by the proof of Lemma~\ref{lem:slice-intersect}. This gives the result. \EPf

\subsection{Reducible representations}

In this section we prove that the copolarity of a direct
sum of representations is not smaller than the 
copolarity of its summand representations. 
The result is:

\begin{thm}\label{thm:red} 
Let $(G,V)$ be an orthogonal representation and suppose that
$V=V_1\oplus V_2$ is an invariant decomposition.
If $\copol{G}{V}\leq k$, then $\copol{G}{V_i}\leq k$ for $i=1$, $2$.
\end{thm}

\begin{lem}\label{lem:p_2}
Given a $(G,V_1)$-regular point $p_1\in V_1$, there is 
$p_2\in V_2$ such that $p=(p_1,p_2)$ is $(G,V)$-regular.
\end{lem}

\Pf Consider the representation $(G_{p_1},V_2)$ and take a
$(G_{p_1},V_2)$-regular point $p_2\in V_2$. We claim that $p=(p_1,p_2)$ is
$(G,V)$-regular. In fact, let $q=(q_1,q_2)\in V_1\oplus V_2$. 
Since $p_1$ is $(G,V_1)$-regular, there is $h\in G$ such that
$G_{p_1}\subset hG_{q_1}h^{-1}=G_{hq_1}$, and therefore 
$(G_{p_1})_{hq_2}\subset(G_{hq_1})_{hq_2}=G_{hq}$. 
Now $p_2$ is $(G_{p_1},V_2)$-regular, so 
there is $g\in G$ such that $G_p=(G_{p_1})_{p_2}\subset 
g(G_{p_1})_{hq_2}g^{-1}\subset gG_{hq}g^{-1}=(gh)G_q(gh)^{-1}$. 
Hence, the result. \EPf

\begin{lem}\label{lem:p_2-Sigma}
Let $\Sigma$ be a $k$-section of $(G,V)$. Given a
$(G,V_1)$-regular point $p_1\in\Sigma\cap V_1$, there is 
$p_2\in V_2$ such that $p=(p_1,p_2)$ is $(G,V)$-regular
and $p\in\Sigma$. 
\end{lem}

\Pf We already know from Lemma~\ref{lem:p_2} that there is $p'_2\in V_2$ 
such that $p'=(p_1,p'_2)$ is $(G,V)$-regular. Let $g\in G$ be 
such that $gp'=(gp_1,gp'_2)\in\Sigma$. Note that $gp_1\in\nu_{gp'}(Gp')$.
Since $gp'$ is $(G,V)$-regular, $\nu_{gp'}(Gp')\subset\Sigma$ 
so that $gp_1\in\Sigma$. Now $p_1\in\Sigma\cap g^{-1}\Sigma$.
Therefore, by Lemma~\ref{lem:trans-sec},
there is $k\in G_{p_1}$ such that $kg^{-1}\Sigma=\Sigma$. 
Hence $p=kp'=(p_1,kp'_2)\in\Sigma$. \EPf

\begin{lem}\label{lem:normal-red}
Let $\Sigma$ be a $k$-section of $(G,V)$.
Given a $(G,V_1)$-regular point $p_1\in\Sigma\cap V_1$,
we have that $\Sigma\cap V_1$ contains $\nu_{p_1}^{V_1}(Gp_1)$
as a subspace of codimension at most $k$,
where $\nu_{p_1}^{V_1}(Gp_1)$ denotes the normal space to $Gp_1$ 
at $p_1$ in $V_1$. 
\end{lem}

\Pf Use Lemma~\ref{lem:p_2-Sigma} to find $p_2\in V_2$
such that $p=(p_1,p_2)$ is $(G,V)$-regular
and $p\in\Sigma$. Now $\Sigma$ contains $\nu_p(Gp)$ with codimension $k$
and $\nu_{p_1}^{V_1}(Gp_1)=\nu_p(Gp)\cap V_1$. \EPf

\medskip

\textit{Proof of Theorem~\ref{thm:red}.}
Fix a $(G,V_1)$-regular $p_1\in V_1$ and fix a $k$-section
$\Sigma$ for $(G,V)$ with $p_1\in\Sigma$. Define 
$\Sigma_1=\cap_{h\in G_{p_1}}h\Sigma\cap V_1$
and let $q_1\in\Sigma_1$ be $(G,V_1)$-regular (for instance,
$q_1$ could be equal to $p_1$). Then, for all
$h\in G_{p_1}$, we have $q_1\in h\Sigma\cap V_1$ and $h\Sigma$ is a 
$k$-section for $(G,V)$. It follows by Lemma~\ref{lem:normal-red}
that $h\Sigma\cap V_1$ contains $\nu_{q_1}^{V_1}(Gq_1)$ as
a subspace of codimension at most $k$ for all $h\in G_{p_1}$.
Therefore $\Sigma_1$ contains $\nu_{q_1}^{V_1}(Gq_1)$ as
a subspace of codimension at most $k$ and hence
$\Sigma_1$ satisfies conditions~(C1), (C2) and~(C3).
This already implies that $\copol{G}{V_1}\leq k$,
but we go on to show that $\Sigma_1$ itself satisfies condition~(C4).

Suppose that
$p_1\in\Sigma_1\cap g^{-1}\Sigma_1$ for some $g\in G$. 
Then $p_1\in h\Sigma\cap g^{-1}h\Sigma$ for all $h\in G_{p_1}$.
We use Lemma~\ref{lem:trans-sec} to find $l=l(h)\in G_{p_1}$ such that 
$lh\Sigma=g^{-1}h\Sigma$. This gives 
$g^{-1}\Sigma_1=\cap_{h\in G_{p_1}}lh\Sigma\cap V_1\supset\Sigma_1$.
Therefore $g\Sigma_1=\Sigma_1$. This shows that 
the distribution $\D_1$ defined by 
$\D_{1q_1}=\Sigma_1\cap T_{q_1}(Gq_1)$ for $q_1\in Gp_1$
satisfies assertion~(c) in Proposition~\ref{prop:distr}
with respect to the principal $(G,V_1)$-orbit $Gp_1$. 
It follows that $\Sigma_1$ satisfies condition~(C4) and hence 
it is a $k_1$-section with $k_1=\dim\D_1\leq\dim\D=k$. \EPf

\subsection{Minimal $k$-sections, the osculating spaces of orbits
and the integrability of $\E$}

In this section we study properties of minimal $k$-sections
of orthogonal representations with nontrivial copolarity
(Recall that an isometric action $(G,M)$ has 
nontrivial copolarity if a minimal $k$-section is a 
proper subset of $M$.)
In particular, we show that in the irreducible
case there can be no minimal $k$-sections
of codimension one or two in the ambient space (Corollary~\ref{cor:cod2}). 
We also characterize the case where the distribution $\E$ is 
integrable (Theorem~\ref{thm:int}). 

Let $(G,V)$ be an orthogonal representation. 
Define an equivalence relation in the set of $G$-regular points by
declaring two points to be equivalent if they can be joined
by a polygonal path which is (at smooth points)
tangent to the distribution of normal 
spaces of the $G$-orbits. 
It is clear that the group action permutes the
equivalence classes. 
For a $G$-regular $p\in V$, denote the equivalence class of $p$ by 
$\mathfrak S_p$. Note that $0\in\mathfrak S_p$, as $p$ is a vector in the 
normal space of $Gp$ at $p$. Therefore the affine hull 
$\vecsp{\mathfrak S_p}$
is a vector subspace of $V$. 

Next suppose that $(G,V)$ admits a $k$-section $\Sigma$.
Let $p\in\Sigma$ be $G$-regular. 
It is clear that $\vecsp{\mathfrak S_p}\subset\Sigma$. 
Let $\Lg$ be the Lie algebra of $G$ and 
consider its induced action by linear skew-symmetric
endomorphisms on $V$. For each 
$X\in\Lg$, let $f_X:\Sigma\to\Sigma$ be defined by 
$f_X(q)=\Pi\circ X_q$, where $\Pi:V\to\Sigma$ is orthogonal projection. 
Then $f_X$ is a linear map and $W_X:=\ker f_X$ is a subspace of $\Sigma$.
Let $p\in\Sigma$ be $G$-regular. 
Define the following subspace of $\Sigma$:
\[ \bar\Sigma_p:=\bigcap_{W_X\ni p}W_X. \] 
We have that $\bar\Sigma_p$ is the subset of $\Sigma$ 
comprised of the common zeros of the vector fields that are
orthogonal projections 
onto $\Sigma$ of the $G$-Killing fields which are orthogonal to $\Sigma$ 
at $p$. It is clear from this definition that if $q\in\bar\Sigma_p$
is a $G$-regular point then $\bar\Sigma_q\subset\bar\Sigma_p$. 

\begin{prop}\label{prop:barsigma}
We have that $\vecsp{\mathfrak S_p}\subset\bar\Sigma_p\subset\Sigma$
and $\bar\Sigma_p$ is an $l$-section of $(G,V)$ with $l\leq k$.
\end{prop}

\Pf Let $X\in\Lg$ be such that $p\in W_X$. Then $X_p\perp\Sigma$. 
If $q\in\mathfrak S_p$, then $q$ can be joined to $p$ by 
a polygonal path which is normal to the orbits. It follows by an iterated
application of Lemmas~\ref{lem:2} and~\ref{lem:3} below that
$X_q\perp\Sigma$, that is, $q\in W_X$. This shows that 
$\mathfrak S_p\subset W_X$. Therefore $\vecsp{\mathfrak S_p}\subset W_X$. 
Since $X$ can be any element in $\Lg$ satisfying $W_X\ni p$, we get 
that $\vecsp{\mathfrak S_p}\subset\bar\Sigma_p$.

We next prove that $\bar\Sigma_p$ is an $l$-section. Since 
$\bar\Sigma_p\subset\Sigma$, it will follow that $l\leq k$. 
In fact, condition (C1) for $\bar\Sigma_p$ is obvious and condition (C2) 
follows from the facts that $\bar\Sigma_p\supset\mathfrak S_p$
and $\mathfrak S_p$ intersects all orbits. Let us verify 
condition (C4). Let $q\in\bar\Sigma_p$ be a $G$-regular point 
and $g\in G$ with $gq\in\bar\Sigma_p$. We need to show that
$g\bar\Sigma_p=\bar\Sigma_p$. 
First we note that it is clear from the definition of $\bar\Sigma_p$ that
$\bar\Sigma_{gq}\subset\bar\Sigma_p$.
Next, since $\mathfrak S_p$ intersects all
orbits and $\mathfrak S_p\subset\bar\Sigma_p$, 
we may assume that $q\in\mathfrak S_p$, and this implies 
as above, via Lemmas~\ref{lem:2} and~\ref{lem:3}, 
that $\bar\Sigma_q=\bar\Sigma_p$. Finally, 
since $q$, $gq\in\bar\Sigma_p\subset\Sigma$, condition 
(C4) for $\Sigma$ gives that $g\Sigma=\Sigma$,
and then we have $gW_X=W_{gXg^{-1}}$, which shows that 
$g\bar\Sigma_q=\bar\Sigma_{gq}$. Putting all this together
we have $g\bar\Sigma_p=g\bar\Sigma_q=\bar\Sigma_{gq}
\subset\bar\Sigma_p$, and hence $g\bar\Sigma_p=\bar\Sigma_p$.

In order to check (C3), let $q\in\bar\Sigma_p$ be a $G$-regular point.
Since $\mathfrak S_p$ intersects all orbits, there is $g\in G$ such that
$gq\in\mathfrak S_p$. It is clear that 
$\nu_{gq}(Gq)\subset\vecsp{\mathfrak S_p}\subset\bar\Sigma_p$. 
Then $\nu_q(Gq)=g^{-1}\nu_{gq}(Gq)\subset g^{-1}\bar\Sigma_p=
\bar\Sigma_p$, where the last equality follows from (C4) 
for $\bar\Sigma_p$. \EPf

\begin{cor}\label{cor:barsigma}
If $\Sigma$ is a \emph{minimal} $k$-section for $(G,V)$ and $p\in\Sigma$
is $G$-regular, then for all $X\in\Lg$ we have that:
$X_p\perp\Sigma$ if and only if $X_q\perp\Sigma$ for all 
$q\in\Sigma$. 
\end{cor}

\Pf This is clear, because $\Sigma=\bar\Sigma_p$. \EPf

\medskip

Fix a minimal $k$-section for $(G,V)$ and a $G$-regular
$p\in\Sigma$. Let $\Lk$ be the Lie subalgebra of $\Lg$
generated by all the $X\in\Lg$ such that $X_p\perp\Sigma$
and let $K$ be the associated connected subgroup of $G$.
Let $G_\Sigma^0$ the connected component of the identity 
in~$G_\Sigma$.
Denote by $G'$ the subgroup of $G$ generated by $G_\Sigma^0$ and 
$K$. Note that the orbits $G'p$ and $Gp$ have the same tangent space
at $p$ (cf.~Lemma~\ref{lem:inj}). 
Therefore they are equal. Since $Gp$ is a principal
orbit, it is not hard to show that $(G',V)$ and $(G,V)$ 
are orbit equivalent (cf.~Lemma~3.6 in~\cite{GTh1}).
By replacing $G$ by $G'$,  
we may now assume that $G$ is generated by $G_\Sigma^0$ and 
$K$. Next we show that $K$ is a normal subgroup of $G$. 
It is enough to show that if $X\in\Lg$ satisfies
$X_p\perp\Sigma$ and $g\in G_\Sigma$, 
then $Y=gXg^{-1}\in\Lk$. In fact,
$Y_{gp}=gX_p\perp g\Sigma=\Sigma$ and hence 
$Y\in\Lk$ by Corollary~\ref{cor:barsigma}.
Since $K$ is normal in $G$, we get that $G$ is a quotient
of the semidirect product 
of $G_\Sigma^0$ and $K$. 

\begin{prop}\label{prop:cod}
Let $(G,V)$ be irreducible with nontrivial copolarity $k$
and assume that $\Sigma$ is a $k$-section. Then, for every
$G$-regular $p\in\Sigma$, there does not exist a nonzero
$\xi\in\nu_p(Gp)$ such that the Weingarten operator 
$A_\xi|_{\E_p}=0$.
\end{prop}

\Pf Suppose there is a nonzero $\xi\in\nu_p(Gp)$
such that the Weingarten operator
$A_\xi|_{\E_p}=0$. Let $\hat \xi$ be the equivariant normal
vector field along $Gp$ which extends $\xi$. Then $\nabla^\perp_u\hat\xi=0$
for all $u\in\E_p$ by Corollary~\ref{cor:partial-parallel} below. 
This implies that $\hat\xi$ 
is constant on $Kp$ as a vector in $V$, so that $\xi$ is in the fix-point
set of $K$. Since $K$ is normal in $G$, the fix-point set of $K$ is 
$G$-invariant. Since $G$ is irreducible on $V$, $K$ must be trivial 
on $V$ and then $\Sigma=V$, but this is impossible as 
$(G,V)$ has nontrivial copolarity. \EPf

\begin{cor}\label{cor:cod1}
Let $(G,V)$ be irreducible with nontrivial copolarity $k$.
Then the cohomogeneity of $(G,V)$ is less than $\frac{l(l+1)}2$, where
$l$ is the codimension of a $k$-section in $V$. 
\end{cor}

\Pf Let $Gp$ be a principal orbit and choose a minimal $k$-section 
$\Sigma\ni p$. 
It follows from Proposition~\ref{prop:cod} that the map 
\[ \xi\in\nu_p(Gp) \mapsto A_\xi|_{\E_p}\in\mbox{Sym}^2(\E_p^*) \]
is injective, where $\mbox{Sym}^2(\E_p^*)$ denotes the symmetric 
square of the dual space of $\E_p$. Now we need just note that
$\dim\E_p=l$. \EPf

\begin{cor}\label{cor:cod2}
Let $(G,V)$ be irreducible with nontrivial copolarity $k>0$
and let $n$ be the dimension of a principal orbit.
Then $k\leq n-3$.
\end{cor}

\Pf If $k=n-1$, then $l=n-k=1$, and by Corollary~\ref{cor:cod1}
the codimension of a principal orbit is $1$. In this case
$G$ is transitive on the unit sphere and therefore polar, so this case is 
impossible.
 
If $k=n-2$, then $l=2$ and the codimension of a principal orbit is at 
most $3$. Here $G$ is either polar on $V$, or has copolarity $1$
and cohomogeneity $3$, by the 
classification of irreducible representations of cohomogeneity 
at most $3$~\cite{Y}. Since $k>0$, we must have $k=1$. Then $n=3$, 
but none of the irreducible representations of cohomogeneity $3$
has principal orbits of dimension $3$. So this case is impossible either. \EPf

\medskip

We next characterize the orthogonal representations 
of nontrivial copolarity whose distribution $\E$ is integrable.

\begin{thm}\label{thm:int}
Let $(G,V)$ be an orthogonal representation with nontrivial
copolarity. Suppose that the distribution $\E$ is integrable. 
Then there is an orthogonal decomposition $V=V_1\oplus V_2$,
a polar representation $(K,V_1)$ and 
another orthogonal representation $(H,V_2)$ such that $(G,V)$ 
is orbit equivalent to 
the direct product representation $(K\times H,V_1\oplus V_2)$.
Here the leaves of $\E$ correspond to the $K$-orbits. 
\end{thm} 

\Pf Let $K$ be the normal subgroup of $G$ as above.
We first show that the leaves of $\E$ coincide with the $K$-orbits. 
For that purpose, note that $\E_q\subset T_q(Kq)$ 
for every $G$-regular $q\in\Sigma$, and then 
$\E_q\subset T_q(Kq)$ for every $G$-regular $q\in V$,
as every $K$-orbit intersects $\Sigma$ and $\E$ is 
$K$-invariant. It follows that,
if $p\in\Sigma$ is a fixed $G$-regular point
and $\beta$ is the leaf of $\E$ through $p$, then $\beta\subset Kp$.
Let $K_\beta^0$ denote the connected component of
the stabilizer of $\beta$ in $K$.
It is clear that $K_\beta^0 p=\beta$. Let $X\in\Lg$ be such that
$X_p=u\in\E_p$. Since $\E_p=T_p\beta=T_p(K_\beta^0 p)$, 
there is an $Y$ in the Lie algebra of $K_\beta^0$
such that $Y_p=u$. Now $Z=X-Y\in\Lk$ and $Z_p=0$. Therefore
the one-parameter subgroup of $K$ generated by 
$Z$ is in the isotropy subgroup $K_p$.
But $K_p\subset K_\beta^0$, since $K$ maps $\E$-leaves
onto $\E$-leaves. It follows that $Z$ is in the Lie algebra 
of $K_\beta^0$ and so is $X$. Since $X$ is an arbitrary generator of 
$\Lk$, it follows that $K_\beta^0=K$ and hence $\beta=Kp$. 

Now it is clear that $\Sigma$ is a section of $(K,V)$ 
so that $(K,V)$ is polar. Note that the $G$-regular points in $V$ are also 
$K$-regular. Let $N$ denote the normalizer of $K$ in the 
orthogonal group $\mathbf O(V)$. Then $N$ maps $K$-orbits onto
$K$-orbits. Let $n\in N$. Since $\Sigma$ intersects all $K$-orbits, there
is $k\in K$ such that $knp\in\Sigma$. Now $kn$ maps $\Sigma$ onto $\Sigma$. 
The principal 
$K$-orbits are isoparametric in $V$ and $kn$ 
preserves their common focal set. 
Therefore $kn$ preserves the focal hyperplanes in $\Sigma$. 
Decompose $\Sigma$ into an orthogonal sum $\Sigma_1\oplus V_2$, 
where $\Sigma_1$ is the span of the curvature normals of the 
principal $K$-orbits. 
If $q\in\Sigma$ is $K$-regular, then $Kq$ is full in the affine subspace 
$q+V_1$, where $V_1$ is the orthogonal complement of 
$V_2$ in $V$, and $K$ acts trivially on $V_2$. 
The decomposition $V=V_1\oplus V_2$ is $N$-invariant. 
Now $kn$ maps $\Sigma_1$ onto $\Sigma_1$ and maps the $K$-orbits in 
$V_1$ onto $K$-orbits in $V_1$. 
Therefore $kn$ is in the Weyl group of $(K,V_1)$, which is 
a finite group generated by the reflections on the focal
hyperplanes in $\Sigma_1$. It follows that $kn$ maps a $K$-orbit
in $V_1$ onto the \emph{same} $K$-orbit in $V_1$. 
In particular, if $n$ is in the connected component $N^0$ of $N$, then 
$kn$ is the identity on $V_1$.
In any case, $N$ and $K$ have the same orbits in $V_1$. 

Since $G$ normalizes $K$, we have
$Kp_1=Gp_1$ for $p_1\in V_1$, and then $G_\Sigma p_1=
Gp_1\cap\Sigma=Kp_1\cap\Sigma_1$ is finite. Now
$G_\Sigma^0p_1=\{p_1\}$. This shows that $G_\Sigma^0$ is trivial 
on $V_1$. Since $K$ is trivial on $V_2$, 
the intersection $G_\Sigma^0\cap K$ is in the kernel 
of the $G$-action on $V$. Therefore we may assume that
$G_\Sigma^0\cap K=\{1\}$. Now the action of $G$ on $V$
is orbit equivalent to the action of the
direct product of $G_\Sigma^0$ and $K$ and this completes the proof of
the theorem. \EPf

\section{Representations of copolarity one}\label{sec:one}
\setcounter{thm}{0}

In this section we describe the structure of a
principal orbit of a representation of copolarity one
(Theorem~\ref{thm:one}) and classify the irreducible 
representations of copolarity one
(Corollary~\ref{cor:one}). 
We start with a lemma with an interest
of its own.

\begin{lem}\label{lem:decomp}
Let $(G,V)$ be an orthogonal representation with  
copolarity $k$. Suppose that 
$p\in V$ is a regular point and $\Sigma$ is a 
$k$-section through $p$. Then, given a Killing field $X$ of $V$ induced by 
$G$, we can write $X=X_1+X_2$, where
$X_1$, $X_2$ are Killing fields induced by $G$ such that
$X_1|_{\Sigma}$ is always tangent to $\Sigma$ and $X_2|_{\Sigma}$ 
is always perpendicular to $\Sigma$.
\end{lem}

\Pf Let $\bar X$ be the intrinsic Killing field of $\Sigma $ which 
is obtained by projecting $X|_{\Sigma}$ to the tangent space of $\Sigma$. 
Since $X_p$ is perpendicular to the normal space $\nu_p(Gp) \subset
\Sigma$, we must have that $\bar X_p$ is tangent to
$G p \cap \Sigma$ at $p$. By Lemma~\ref{lem:inj}, 
$G_\Sigma $ acts transitively on $Gp \cap \Sigma$, 
so there exists a Killing field $X_1$ induced by
$G_\Sigma \subset G$ such that $X_{1p} = \bar X_p$ (observe that the 
restriction to $\Sigma$ of $X_1$ is always tangent to $\Sigma$). 
Then $X_2 = X-X_1$ is perpendicular to $\Sigma $ at $p$, and, 
by Corollary~\ref{cor:barsigma}, we 
have that $X_2|_{\Sigma}$ is always perpendicular to $\Sigma$. \EPf

\begin{thm}\label{thm:one}
Let $(G,V)$ be an orthogonal representation with copolarity 
$k=1$ and let $M = Gp$ be a principal orbit. Then the submanifold $M$ of 
$V$ splits extrinsically as $M = M_0 \times M_1$,
where $M_0$ is a homogeneous isoparametric submanifold and $M_1$ is
either one of the following:
\begin{enumerate}
\item[(i)] a nonisoparametric homogeneous curve;
\item[(ii)] a focal manifold of an irreducible
homogeneous isoparametric 
submanifold which is obtained by focalizing a one-dimensional 
distribution;
\item[(iii)] a codimension $3$ homogeneous 
submanifold.  
\end{enumerate}
\end{thm}

\Pf It follows from the proof of Theorem~B in~\cite{OS}
that the Lie algebra $\Lh$, which is
(algebraically) generated by the projection to 
$\nu_p(M)$ of Killing fields induced by $G$ (restricted to $\nu_p (M)$), 
contains the normal holonomy algebra (and it is contained in its normalizer). 
By Lemma~\ref{lem:decomp}, 
$\Lh$ is generated by the projection of the 
Killing fields induced by $G_\Sigma$, where $\Sigma$ is a $1$-section 
through $p$. But if $X$, $Y\neq 0 $ are such Killing fields then they  
must be proportional at $p$, since $\dim (Gp\cap\Sigma)=1$. By multiplying 
one of them by a nonzero scalar, we may assume that $X_p=Y_p$. 
Now, by Corollary~\ref{cor:barsigma},
$(X-Y)|_{\Sigma}$ must be always perpendicular 
to $\Sigma$ and is thus zero. 
Therefore $\Lh$ is generated by $X|_{\Sigma}$ and so it 
has dimension $1$. So, the restricted normal holonomy group of $M$ has 
dimension $0$ or $1$. 

Let $\bar M$ be a nonisoparametric irreducible extrinsic factor of $M$.
We may assume it is full 
(since we are only concerned with the geometry of $M$).
If $\bar M$ has flat normal bundle, then, by~Theorem~A in~\cite{O1},
$\bar M$ is a homogeneous curve (which is not an extrinsic circle).
Assume that $\bar M$ has
nonflat normal bundle. 
Orthogonally decompose the normal bundle
\[ \nu (\bar M) = \nu _0(\bar M)  \oplus \nu _s(\bar M), \]
where $\nu _0(\bar M)$ is the maximal parallel and flat subbundle of
$\nu (\bar M)$. By the Normal Holonomy Theorem~\cite{O2}, as
the restricted normal holonomy group of $\bar M$ 
has dimension $1$, $\nu_s(\bar M)$ has dimension $2$ over $\bar M$ 
(and the restricted normal 
holonomy group acts as the circle action in a two-dimensional Euclidean 
space); notice that there can be at most one irreducible 
extrinsic factor 
of $M$ with nonflat normal bundle. 
If the codimension of $\bar M$ is greater than $3$, then  
$\mbox{rank}(\bar M)$  
(i.e. the dimension over $\bar M$ of $\nu _0(\bar M))$ is 
at least $2$. Then, by Theorem~A in~\cite{O1}, 
$G$ can be enlarged to a group admitting a representation which is
the isotropy representation of an
irreducible symmetric space and has $\bar M$ as an orbit. 
The normal holonomy tube, 
which has one dimension more and coincides with a principal orbit of the 
isotropy representation, is an irreducible isoparametric submanifold.
We finish the proof by observing that there cannot be 
two distinct nonisoparametric irreducible extrinsic factors of $M$,
for otherwise the copolarity would be at least $2$. 
\EPf

\begin{cor}\label{cor:one}
Let $(G,V)$ be an irreducible representation 
of nontrivial copolarity $1$. Then $(G,V)$ is 
one of the three orthogonal representations
listed in the table of Theorem~\ref{thm:1}.
\end{cor}

\Pf 
We know from 
Theorem~\ref{thm:one} that any 
principal orbit is either 
a codimension $3$ homogeneous 
submanifold or a focal 
manifold of an irreducible
homogeneous isoparametric submanifold.
If there is a principal orbit which falls into the first case,
then the cohomogeneity is $3$ and the result follows from the classification
of cohomogeneity $3$ irreducible representations (see~\cite{Y}). 
Suppose, on the contrary, that no principal orbit
falls into the first case. 
Then any principal orbit is a focal 
manifold of an isoparametric submanifold,
thus it is taut (see~\cite{HPT}).
If we can show that the nonprincipal orbits are also taut, then 
it will follow from the classification of taut irreducible 
representations~\cite{GTh3}
that the cohomogeneity is $3$ and this is a contradiction. 


Let $Gq$ be a 
nonprincipal orbit. We need to prove that it
is tautly embedded in $V$. 
We can find a vector $v\in\nu_q(Gq)$
and a decreasing sequence $\{t_n\}$ such that $t_n\to0$ and 
$p_n=q+t_nv$ are $G$-regular points. For each $n$ there is a
group $K_n\supset G$ and an isotropy representation 
of a symmetric space
$(K_n,V)$ such that $Gp_n=K_np_n$. Since the number of 
isotropy representations of symmetric spaces of a given
dimension is finite, by passing to a subsequence 
we can assume that there is a sequence $h_n\in\SO{V}$ 
such that $K_n=h_nK_1h_n^{-1}$ for $n\geq2$. By compactness of $\SO V$,
again by passing to a subsequence we can write $h_n\to h\in\SO V$. 
Let $K_\infty=hK_1h^{-1}$. We finally prove that $K_\infty q=Gq$.
Since $(K_\infty, V)$ is conjugate to the 
isotropy representation of a symmetric space,
its orbits are tautly embedded \cite{BS}, and this 
shows that $Gq$ is taut. 

Let $k\in K_\infty$. We have that $k=hk_1h^{-1}$ for some
$k_1\in K_1$, and then $k=\lim k_n$, where $k_n=h_nk_1h_n^{-1}\in K_n$
for $n\geq2$.
Since $K_np_n=Gp_n$, there is $g_n\in G$ such that $g_np_n=k_np_n$.
By passing to a subsequence we may assume that $g_n\to g\in G$.
Now $kq=\lim k_np_n=\lim g_np_n=gq\in Gq$, and this proves that
$K_\infty q\subset Gq$. Since the reverse inclusion
is clear, this completes the proof of the claim and the proof
of the corollary. \EPf

\begin{rmk}
\em
Regarding Theorem~\ref{thm:one} and
Corollary~\ref{cor:one}, we know examples
of reducible representations of nontrivial
copolarity $1$, but all of them have 
cohomogeneity $3$. Thus we do not know if case 
(ii) in Theorem~\ref{thm:one} indeed can occur. 
\end{rmk}

\section{Variational co-completeness}\label{sec:varcomp}
\setcounter{thm}{0}

Let $N$ be a submanifold of a complete Riemannian manifold~$M$. 
Let $\eta:\nu(N)\to M$ denote the \emph{endpoint map} of $N$, that is,
the restriction of the exponential map of $M$ to the normal bundle of 
$N$. 
A point $q=\eta(v)$ is a \emph{focal point of $N$ in the direction 
of $v\in\nu(N)$of multiplicity $m>0$} if $d\eta_v:T_v\nu(N)\to T_qM$ is
not injective and the dimension of its kernel is $m$.
Let $v\in\nu_p(N)$ and $\gamma_v$ denote the geodesic $t\mapsto\exp_p(tv)$.
A Jacobi field along $\gamma_v$ is called an $N$-\emph{Jacobi field}
if it is the variational vector field of a variation through geodesics
that are at time zero orthogonal to $N$. 
We will denote the space of $N$-Jacobi 
fields along $\gamma_v$ by $\mathcal{J}^N(\gamma_v)$. 
It is not difficult to see that $J$ is an
$N$-Jacobi field along $\gamma_v$ if and only if $J(0)\in T_pN$ and
$J'(0)+A_vu\in\nu_p(N)$, where $p$ is the footpoint of $v$, $u=J(0)$ 
and $A_v$ is the Weingarten map in direction $v$. 
The point $q$ is a focal
point of $N$ in the direction $v$ if there is an $N$-Jacobi 
field along $\gamma_v$ that vanishes in $q$. 
We will denote the space of $N$-Jacobi 
fields along $\gamma_v$ that vanish in $q$
by $\mathcal{J}^N_q(\gamma_v)$. 

Now let $(G,M)$ be an isometric action of a compact Lie group 
$G$ on the complete Riemannian manifold $M$.
The action $(G,M)$ is called \emph{variationally complete}
if every Jacobi field $J\in\mathcal{J}^N_q(\gamma_v)$, 
where $N$ is a $G$-orbit and $q$ is a focal point
of $N$ in the direction of $v$, is the restriction along $\gamma_v$ of 
a Killing field on $M$ induced by the action of $G$.

More generally, let $N$ be a fixed principal orbit of an isometric 
action $(G,M)$ and let $p\in N$. For each $v\in\nu_p(N)$,
we have an isomorphism $\mathcal{J}^N(\gamma_v)\to
T_p N\oplus\nu_p N=T_p M$ given by $J\mapsto(J(0),J'(0)+A_vJ(0))$.
Let $U_p$ be a subspace of $T_p M$
with the following property:
\begin{enumerate}
\item[(P)] for each $v\in\nu_p(N)$,
if $J\in\mathcal{J}^N(\gamma_v)$ vanishes for some
$t_0>0$ and $(J(0),J'(0)+A_vJ(0))$ is orthogonal to 
$U_p$, then $J$ is the restriction along $\gamma_v$ of a 
$G$-Killing field.
\end{enumerate}
Of course, $U_{gp}$ can always be taken to be 
$g_*U_p$, where $g\in G$, and in particular, $U_p$
can always be taken to be $G_p$-invariant. 
Moreover, $U_p$ can always be taken to be
all of~$T_p M$.  
In any case we write $\mbox{covar}_N(G,M)\leq\dim U_p$.
We say that the \emph{variational co-completeness} of $(G,M)$
is less than or equal to $k$, where $k$ is an integer between $0$ and 
$\dim M$, and we write $\covar{G}{M}\leq k$,
if $\mbox{covar}_N(G,M)\leq k$ for all principal $G$-orbits $N$. 
Clearly, 
$\covar{G}{M}\leq0$ (or, in this case, $\covar{G}{M}=0$) 
if and only if $(G,M)$ is 
variationally complete, and one cannot do better
than $\covar{G}{M}\leq\dim M$
for a generic isometric action. 
In the next section we will describe a situation where the 
intermediate values occur. 

Observe that the intersection of two subspaces of $T_pM$
with property~(P) does not need to have that property,
so we cannot speak of a minimal subspace with property~(P). 
Nonetheless, given a general isometric action $(G,M)$ and a
$G$-regular $p\in M$, 
we next show how to construct a canonical subspace $U_p^0$ of $T_p M$ with 
property~(P). Let $N=Gp$. 
For each $v\in\nu_p(Gp)$ and
$q$ a focal point of $N$ in the direction $v$,
consider the subspace 
$\tilde U_p^{v,q}$ of $T_pM$ spanned by the initial conditions
$(J(0),J'(0)+A_vJ(0))$ for all $J\in\mathcal{J}^N_q(\gamma_v)$.
Now take the subspace of $\tilde U_p^{v,q}$ spanned 
by the initial conditions
of all $G$-Killing fields in $\mathcal{J}^N_q(\gamma_v)$
and let $U_p^{v,q}$ denote its orthogonal complement
in~$\tilde U_p^{v,q}$.
Finally, define $U_p^0$ as the sum over 
$v$, $q$ of the subspaces $U_p^{v,q}$. 
It is clear that $U_p^0$ has property~(P). 
Note also that these $U_p^0$, for $p\in N$, define a $G$-invariant 
distribution on $N$.

\subsection{The theorem of Conlon for actions admitting $k$-sections}

In this section we prove a version of Conlon's
theorem~\cite{C}.

\begin{thm}\label{thm:Conlon}
If $(G,M)$ is an isometric action 
admitting a $k$-section which is flat in the induced 
metric, then $\covar{G}{M}\leq k$.
\end{thm}

Let $N=Gp$ be a principal orbit, $v\in\nu_p(N)$,
and choose a flat $k$-section $\Sigma$ through $p$. 

\begin{lem}\label{lem:1}
Let $J$ be an $N$-Jacobi field along $\gamma_v$ such that $J(0)\in\E_p$.
If $J(t_0)=0$ for some $t_0>0$, then $J$ is always orthogonal to $\Sigma$.
\end{lem}

\Pf Decompose $J=J_1+J_2$, where $J_1(t)$ and $J_2(t)$ are respectively 
the tangent and normal components of $J(t)$ relative to 
$T_{\gamma_v(t)}\Sigma$. Since $\Sigma$ is totally geodesic, 
$J_1$ and $J_2$ are Jacobi fields along $\gamma_v$. 
Now $J_1$ is a Jacobi field in $\Sigma$ with
$J_1(0)=J_1(t_0)=0$. Since $\Sigma$ is flat, we have that $J_1$ identically 
vanishes and hence $J=J_2$. \EPf

\begin{lem}\label{lem:2}
Let $J$ be an $N$-Jacobi field along $\gamma_v$ such that $J(0)\in\E_p$.
If $J$ is the restriction along $\gamma_v$ of a $G$-Killing field on~$M$,
then $J$ satisfies $J'(0)+A_vJ(0)=0$.
\end{lem}

\Pf Let $X$ be a $G$-Killing field on $M$ which restricts to $J$ 
along~$\gamma_v$. Note that $X_p=J(0)\in\E_p$.
Denote by $\nabla$ the Levi-Civita
connection of $M$. Now 
$J'(0)=(\nabla_v X)_p$. 
Let $E(t)$ be any vector field along $\gamma_v(t)$
which is normal to $G\gamma_v(t)$. 
Then $E(t)$ is tangent to $\Sigma$ for small $t$. We have 
$\inn{J'(0)}{E}_p=\inn{\nabla_v X}{E}_p=-\inn{X}{\nabla_v E}_p=0$. 
Therefore $J'(0)\in T_pN$, and this implies $J'(0)+A_vJ(0)\in T_pN$.
Since we already have $J'(0)+A_vJ(0)\in\nu_pN$, we get that
$J'(0)+A_vJ(0)=0$. \EPf

\begin{lem}\label{lem:3}
Let $J$ be an $N$-Jacobi field along $\gamma_v$ such that $J(0)\in\E_p$.
Then $J$ is always orthogonal to $\Sigma$ if and only if
$J$ satisfies $J'(0)+A_vJ(0)=0$.
\end{lem}

\Pf If $J$ is always orthogonal to $\Sigma$, then 
$J'(0)+A_vJ(0)$ is also orthogonal to $\Sigma$ as 
$\Sigma$ is totally geodesic. But 
as an $N$-Jacobi field, $J$ satisfies $J'(0)+A_vJ(0)\in\nu_p(N)$. 
Now we have $\nu_p(N)\subset T_p\Sigma$, hence $J'(0)+A_vJ(0)=0$.
Conversely, if $J'(0)+A_vJ(0)=0$, then $J'(0)=-A_vJ(0)\in\E_p$,
since $\E_p$ is $A_v$-invariant and $J(0)\in\E_p$. Since $J$ and $J'$ are 
both orthogonal to $\Sigma$ at time zero and $\Sigma$ is totally geodesic, 
we deduce that $J$ is always orthogonal to $\Sigma$. \EPf

\medskip

We finish the proof of Theorem~\ref{thm:Conlon} by observing that
$\D_p$ has property (P). 
In fact, if an $N$-Jacobi field $J$ along 
$\gamma_v$ for some $v\in\nu_p N$ 
vanishes for some $t_0>0$ and $(J(0),J'(0)+A_vJ(0))$ 
is orthogonal to $\D_p$, then $J(0)\in\E_p$. By Lemmas~\ref{lem:1}
and~\ref{lem:3}, $J'(0)+A_vJ(0)=0$. Let $X$ be a $G$-Killing on $M$
such that $X_p=J(0)$. Then $X$ restricts to a Jacobi field along
$\gamma_v$ which by Lemma~\ref{lem:2}
must be~$J$. 

\medskip

Since Lemmas~\ref{lem:1}, \ref{lem:2} and~\ref{lem:3} do not depend 
on condition~(C4) in the definition of a $k$-section, we
have the following corollary of the proof. 

\begin{cor}\label{cor:partial-vc}
Let $(G,M)$ be an isometric action.
Suppose there is a flat, connected, complete submanifold $\Sigma$ of $M$ 
satisfying conditions~(C1), (C2) and~(C3) in the definition
of $k$-section. Let $N$ be a principal orbit, $p\in N\cap\Sigma$
and $v\in\nu_pN$. Then $T_p(N\cap\Sigma)$ has property (P).
\end{cor}

Let $N$ be a principal orbit of
an isometric action $(G,M)$ admitting a $k$-section $\Sigma$.
Let $\xi$ be a normal vector field parallel along a curve in $N$
that is everywhere tangent to the 
distribution $\E$. The next corollary implies 
that the principal curvatures of $N$ along $\xi$  
are constant.

\begin{cor}\label{cor:partial-parallel}
Let $(G,M)$ be an isometric action admitting a $k$-section $\Sigma$
(non necessarily flat). Let $N$ be a principal orbit, $p\in N\cap\Sigma$
and $v\in\nu_pN$. Extend $v$ to an equivariant normal vector field 
$\hat v$ along $N$. Then $\hat v$ is parallel along $\E$. 
\end{cor}

\Pf By homogeneity, it is enough to show that $\nabla^\perp_u\hat v=0$ 
for all $u\in\E_p$. Let $X$ be a $G$-Killing field such that 
$X_p=u\in\E_p$. 
Let $J$ be the $N$-Jacobi field along $\gamma_v$ 
which is the restriction of $X$. Then $J(0)=u\in\E_p$.
By Lemma~\ref{lem:2}, we have $J'(0)+A_vu=0$. On the other hand, 
since $(L_X\hat v)_p=0$,
$J'(0)=(\nabla_v X)_p=
(\nabla_{X}\hat v)_p=-A_vu+\nabla_u^\perp\hat v$. 
Hence, $\nabla^\perp_u\hat v=0$. \mbox{}\hfill\EPf

\subsection{A weak converse for~5.1 in the Euclidean case}

In this section we prove a sort of converse to Theorem~\ref{thm:Conlon}
in the Euclidean case (Theorem~\ref{thm:converse})
and obtain, as a corollary, a generalization of a result
about polar representations (Corollary~\ref{cor:normal}).
First, observe that if $(G,V)$ 
is an orthogonal representation and $N=Gp$ is a principal 
orbit, then an~$N$-Jacobi field~$J$
along~$\gamma_v$ which vanishes for some~$t_0>0$ is necessarily 
of the form $J(t)=(1-\frac{t}{t_0})J(0)$ (since the Jacobi equation
in Euclidean space is $J''=0$). Therefore the vector
$J'(0)+A_vJ(0)=(A_v-\frac{1}{t_0}\mbox{id}_{T_p N})J(0)$ is 
simultaneously normal and tangent to $N$,
so that it must vanish. In this way we see that
an element~$J$ of $\mathcal{J}^N_q(\gamma_v)$ is
completely determined by the value of $J(0)$,
so that it is enough to consider property (P)
for subspaces of $T_p N$. 
Note also that the span in $T_pN$ of 
the $J(0)$, where $J\in\mathcal{J}^N_q(\gamma_v)$ and
$J=X\circ\gamma_v$ for some $X\in\Lg$,
is $T_p(G_qp)$. We conclude that in the Euclidean
case property (P) can be rewritten in the following form:
\begin{enumerate}
\item[($\mbox{P}_{\mbox{\scriptsize Euc}}$)] for each $v\in\nu_p(N)$,
if $A_vu=\lambda u$ for some $\lambda\neq0$
and $u$ is orthogonal to $U_p$, then $u\in T_p(G_qp)$,
where $q=p+\frac{1}{\lambda}v$.
\end{enumerate}
Moreover, it follows that
the canonical subspace $U_p^0$ is the sum, 
over $v\in\nu_pN$ and $q$ a focal point of $N$ in the direction $v$, of
the orthogonal complement of $T_p(G_qp)$ in 
$\{J(0):J\in\mathcal{J}^N_q(\gamma_v)\}$.
Recall that $U_p^0$ is $G_p$-invariant, but it is not clear that
$U_p^0$ is invariant under the second fundamental form of $N$. 

\begin{prop}\label{prop:can}
Every subspace $\D_p$ of $T_p N$
which satisfies property ($\mbox{P}_{\mbox{\scriptsize Euc}}$)
and is invariant under the second fundamental form of $N$
must contain $U_p^0$.
\end{prop}

\Pf Suppose, on the contrary, that $U_p^0\not\subset\D_p$.
Then, by the definition of $U_p^0$, there exists 
$v\in\nu_pN$ and an eigenvector
$u$ of $A_v$ with eigenvalue $\lambda\neq0$ such that 
$u\not\in\D_p$ and 
$u$ is in the orthogonal complement of $T_p(G_qp)$ in
$\{J(0):J\in\mathcal{J}^N_q(\gamma_v)\}$, where 
$q=p+\frac{1}{\lambda}v$.

Write $u=u_1+u_2$, where $u_1\in\D_p$ and $u_2\perp\D_p$.
Then $u_2\neq0$. Since $\D_p$ is $A_v$-invariant,
we get that $A_v u_2=\lambda u_2$. By 
($\mbox{P}_{\mbox{\scriptsize Euc}}$) for $\D_p$, we have that
$u_2\in T_q(G_q p)$. Therefore $u_2\perp u$, but this is a 
contradiction to the fact that $u_2\neq0$. \EPf

\begin{thm}\label{thm:converse}
Let $(G,V)$ be an orthogonal representation, $N$ be
a principal orbit and suppose there is a $G$-invariant, 
$k$-dimensional distribution $\D$ in $N$ which is 
autoparallel in $N$ and
invariant under the second fundamental form of $N$,
and satisfies property ($\mbox{P}_{\mbox{\scriptsize Euc}}$).
Then $\copol{G}{M}\leq k$ (and hence, by Theorem~\ref{thm:Conlon},
we have that $\covar{G}{M}\leq k$).
\end{thm}

\Pf Fix $p\in N$. We will prove that $\Sigma=\D_p\oplus\nu_p(Gp)$  
satisfies conditions (C1), (C2) and (C3) for $(G,V)$. 
For that purpose, 
we will use Proposition~\ref{prop:suf}(a). Let $v\in\nu_p N$
be such that the principal curvatures of $A_v$ are all
nonzero (note that the subset of all such $v$
in $\nu_p N$ is open and dense).
Suppose that $q=p+v$ is a $G$-regular point. 
Let $\E$ be the orthogonal complement distribution of $\D$
in $N$. Let $u\in\E_p$.
Since $\D$ is $A_v$-invariant, we may assume that 
$u$ is an eigenvector of $A_v$, and we know that
the corresponding eigenvalue $\lambda$ is not zero.  
Now the $N$-Jacobi field $J$ along $\gamma_v(t)=p+tv$ 
with initial conditions 
$J(0)=u$, $J'(0)+A_v u=0$ is given by $J(t)=(1-t\lambda)u$. 
By property ($\mbox{P}_{\mbox{\scriptsize Euc}}$), $J$ is the restriction of a $G$-Killing
field along~$\gamma_v$. In particular, $J(1)$ is tangent to $T_q(Gq)$. 
Since $Gq$ is a principal orbit, the slice representation at $q$
is trivial, so the one-parameter subgroup of $G$ that induces $J$ 
cannot fix $q$ and thus $J(1)\neq0$. This implies that $u\in T_q(Gq)$. 
We have shown that $\E_p\subset T_q(Gq)$
in the case $q=p+v$ is a $G$-regular point and 
$v\in\nu_p N$ is such that the principal curvatures of $A_v$ are all
nonzero. The case of an arbitrary $v\in\nu_pN$ follows from 
a limiting argument. 
This implies that 
$\Sigma$ satisfies conditions~(C1), (C2) and~(C3)
by Proposition~\ref{prop:suf}(a). \EPf

\begin{rmk}
\em
In Theorem~\ref{thm:converse}, if it happens that
the distribution $\D$ coincides with the distribution defined
by the canonical subspaces, namely $\D_p=U_p^0$ for all $p\in N$, then
we can show that $\Sigma$ satisfies condition (C4) so that
it is already a $k$-section (not necessarily minimal). 
In fact, following the notation of the proof,
suppose that $v\in\nu_pN$ is arbitrary and $q=p+v$ is a $G$-regular point
where $q=gp$ for some $g\in G$. We want to show that $\E_p=\E_q$.
Note that 
$T_q(\Sigma\cap N)$ is invariant under the second fundamental form of 
$N$ at $q$, since $\Sigma$ is totally geodesic.
Moreover, $T_q(\Sigma\cap N)$ has property 
($\mbox{P}_{\mbox{\scriptsize Euc}}$) by 
Corollary~\ref{cor:partial-vc}. Now 
Proposition~\ref{prop:can}
implies that $T_q(\Sigma\cap N)\supset\D_q$.
Since $\E_p$ is the orthogonal complement of $T_q(\Sigma\cap N)$
in $T_qN$, we deduce that $\E_p\subset\E_q$, and thus,
by dimensional reasons, $\E_p=\E_q$.
It follows from 
Proposition~\ref{prop:suf}(b) that 
$\Sigma$ is a $k$-section. 
\end{rmk}

It is known that if a principal orbit of an
orthogonal representation has the 
property that equivariant normal vector fields are parallel
in the normal connection, then the representation is polar. 
The following corollary is a generalization of this 
result.

\begin{cor}\label{cor:normal}
Let $(G,V)$ be an orthogonal representation, $N$ be
a principal orbit and suppose there is a $G$-invariant, 
$k$-dimensional distribution $\D$ in $N$ which is 
autoparallel in $N$ and 
invariant under the second fundamental form of $N$.
Let $\E$ be the orthogonal complement distribution of $\D$
in $N$.
Assume that every equivariant normal vector field
on $N$ is parallel along $\E$.
Then $\copol{G}{M}\leq k$.
\end{cor}

\Pf It is enough to see that $\D$ has property 
($\mbox{P}_{\mbox{\scriptsize Euc}}$) and use Theorem~\ref{thm:converse}.
Let $p\in N$, $v\in\nu_p(N)$ and suppose that
$A_vu=\lambda u$ where $\lambda\neq0$ and
$u$ is orthogonal to $\D_p$. Set $q=p+\frac{1}{\lambda}v$.
Consider the equivariant 
normal vector field $\hat v$ that extends $v$. 
Then $\nabla^{\perp}_u\hat v=0$ by hypothesis. 
Let $X\in\Lg$ be such that $X_p=u$. 
Now \[X_q=\frac{d}{ds}\Big|_{s=0}(\exp sX)q=
\frac{d}{ds}\Big|_{s=0}(\exp sX)p+
\frac{1}{\lambda}\frac{d}{ds}\Big|_{s=0}\hat v(s)
=u+\frac{1}{\lambda}(-A_v u+\nabla^{\perp}_u\hat v)=0.\]
Hence $u\in T_q(G_qp)$. \EPf

\bigskip
\textbf{Final questions}
\begin{enumerate}
\item Is it true that a minimal $k$-section $\Sigma$ of an orthogonal
representation $(G,V)$, where $G$ is the maximal 
(not necessarily connected) subgroup of $\mathbf O(V)$ with its orbits,  
always coincides with the fixed point set of 
a principal isotropy group at a point $p\in\Sigma$?
\item Is there an example of a focal manifold of an
irreducible homogeneous 
isoparametric submanifold, obtained by focalizing a one-dimensional 
distribution, which is an extrinsic factor of a principal orbit 
of a \emph{reducible}
representation of nontrivial
copolarity $1$ and cohomogeneity bigger than $3$?
\item Classify representations of 
nontrivial copolarity (in the irreducible case we believe
that there should be not too many examples). 
\end{enumerate}

\bibliographystyle{amsplain}
\bibliography{paper}

\bigskip

\parbox[t]{7cm}{\footnotesize\sc Instituto de Matem\'atica e Estat\'\i
stica\\
                Universidade de S\~ao Paulo\\
                Rua do Mat\~ao, 1010\\
                S\~ao Paulo, SP 05508-090\\
                Brazil\\ 
                E-mail: {\tt gorodski@ime.usp.br}}

\bigskip

\parbox[t]{9cm}{\footnotesize\sc Facultad de Matem\'atica, Astronom\'\i a
y F\'\i sica\\
Universidad Nacional C\'ordoba\\
Medina Allende y Haya de la Torre\\ 
Ciudad Universitaria\\
5000 C\'ordoba Argentina\\ 
E-mail: {\tt olmos@mate.uncor.edu}}

\bigskip

\parbox[t]{7cm}{\footnotesize\sc Departamento de Matem\'atica\\  
                                 Universidade Federal de S\~ao Carlos\\ 
                                 Rodovia Washington Luiz, km 235\\ 
                                 S\~ao Carlos, SP 13565-905\\
                                 Brazil\\ 
                                 E-mail: {\tt tojeiro@dm.ufscar.br}}

\end{document}